\documentclass[reqno,12pt]{amsart}
\usepackage{amsmath, amsfonts, amsthm, amssymb, verbatim,dsfont}
\usepackage[mathcal]{euscript}
\newtheorem{definition}{Definition}[section]
\newtheorem{lemma}[definition]{Lemma}
\newtheorem{theorem}[definition]{Theorem}

\newtheorem{corollary}[definition]{Corollary}

\numberwithin{equation}{section}
\newcommand \be         {\begin{equation}}
\newcommand \ee         {\end{equation}}
\newcommand \RR     {\mathbb{R}}
\newcommand \Rn     {\mathbb{R}^n}
\newcommand \eps    {\epsilon}
\newcommand \Ocal       {\mathcal O}
\newcommand \Ccal       {\mathcal C}
\newcommand \Scal       {\mathcal S}
\newcommand \Ncal   {\mathcal N}
\newcommand \xbar   {\overline{x}}
\newcommand \Deltabar   {\overline{\Delta}}
\newcommand \del        {{\partial}}
\newcommand \MM     {{\mathbf M^n}}
\newcommand \Lorentz   {{\mathbf M^{n+1}}}
\newcommand \Lorentzplus   {{\mathbf M_+^{n+1}}}
\newcommand \ra         {\big\rangle}
\newcommand \la         {\big\langle}
\newcommand \sgn    {\mbox{sgn}}
\newcommand \Fb         {\widetilde F}
\newcommand \Ub     {\widetilde U}
\newcommand \Hcal   {\mathcal H}
\newcommand \dive   {\mbox{div}_g\hskip.06cm}

\newcommand \grade  {\mbox{grad}_g\hskip.06cm}
\newcommand \ubar   {{\bar u}}
\newcommand \vbar   {{\bar v}}
\newcommand \loc        {\text{loc}}
\begin{document}
\bibliographystyle{plain}
\title[Geometry compatible conservation laws on  manifolds]
{\small Well-posedness theory
\\
for geometry compatible hyperbolic conservation laws on manifolds
}
\author[M. Ben-Artzi and P.G. LeFloch]{Matania Ben-Artzi$^1$ \and Philippe G. LeFloch$^2$
}
\thanks{$^1$ Institute of Mathematics, Hebrew University, Jerusalem 91904, Israel.
E-mail: mbartzi@math.huji.ac.il.
\newline
$^2$ Laboratoire Jacques-Louis Lions \& CNRS UMR 7598, ?
University of Paris VI, 75252 Paris, France. E-mail : lefloch@ann.jussieu.fr.}
\thanks{ 2000\textit{\ AMS Subject Classification:} 35L65,58J45,76N10.
}
\date{\today. \quad To appear in Ann. Inst. H. Poincar\'e, Nonlin. Anal.}
\begin{abstract}
Motivated by many applications (geophysical flows, general relativity), we attempt to set the foundations
for a study of entropy solutions to nonlinear hyperbolic conservation laws posed on a
(Riemannian or Lorentzian) manifold.
The flux of the conservation laws is viewed as a vector-field on the manifold and
depends on the unknown function as a parameter.
We introduce notions of entropy solutions in the class of bounded measurable functions
and in the class of measure-valued mappings.
We establish the well-posedness theory for conservation laws on a manifold, by generalizing
both Kruzkov's and DiPerna's theories originally developed in the Euclidian setting.
The class of {\sl geometry-compatible} (as we call it) conservation laws is singled out
as an important case of interest, which leads to robust $L^p$ estimates
independent of the geometry of the manifold. On the other hand,
general conservation laws solely enjoy the $L^1$ contraction property and leads to a unique contractive
semi-group of entropy solutions. Our framework allows us
to construct entropy solutions on a manifold via the vanishing diffusion method or the finite volume method.
\end{abstract}
\maketitle

\section{Introduction}
\label{1-0}

The theoretical work on discontinuous solutions to nonlinear hyperbolic conservation laws has
been restricted so far to problems set in the Euclidian space $\Rn$. Motivated by
numerous applications, to geophysical fluid flows (shallow-water equations on the surface on the Earth)
and general relativity (Euler-Einstein equations describing neutron stars and black holes) in particular,
we attempt in the present
paper to set the foundations for a study of weak solutions (including shock waves) to
hyperbolic conservation laws posed on a Riemannian or Lorentzian manifold. Recall that
Kruzkov's theory \cite{Kruzkov} deals with equations posed in the Euclidian space $\Rn$
and provides existence, uniqueness, and stability of $L^1 \cap L^\infty$
entropy solutions to nonlinear hyperbolic balance laws in
several space dimensions. In addition, in the case of conservation laws with
``constant flux'', depending solely on
the conservative variable but not on the time and space variables, solutions are known to satisfy
the maximum principle as well as
the $L^1$ contraction and total variation diminishing (TVD) properties.

The above stability features play a central role in the theory and numerical
analysis of hyperbolic conservation laws. The purpose of this paper is to introduce a suitable
generalization of Kruzkov's theory when solutions are defined on a manifold rather than on $\Rn$.
An investigation of the interplay between the geometry of manifolds with limited regularity
and the analysis of hyperbolic partial differential equations (for non-smooth functions)
appears to be particularly challenging.
We hope that the framework initiated here will lead to further works on this important subject.

Let $(\MM,g)$ be an oriented, compact, $n$-dimensional Riemannian manifold.
One of our tasks will be to define a suitable class of conservation laws posed on $\MM$
to which the stability properties of the Euclidian space could possibly extend.
It should be observed that, on a manifold, there is no concept
of ``constant flux'' --in the sense that the flux would be ``independent'' of the space
(and time) variable, since as we will see below the flux is a section of the tangent bundle $T\MM$.
Indeed, the flux at a point $x \in \MM$ belongs to the tangent space $T_x\MM$ to the manifold
at that point.  The Euclidian space is special in that
it is possible to choose the flux $f(u)$ at an arbitrary point $x \in \Rn$ and to parallel-transport it
to the whole space $\Rn$; the vector field generated in this fashion remains constant in $x$
and, therefore, smooth.
By contrast, if a tangent vector is selected at one point of the manifold one can not, in general, associate with it
by parallel-transport a unique smooth vector field on the manifold; this is a consequence of
the curvature of $\MM$. As an example, we recall that any vector field on $S^2$
(the unit sphere in $\mathbb{R}^3$, an important example arising in geophysical flows and
motivating our study) must have two critical points, at least, or else must be discontinuous.

An outline of this paper is as follows.
In Section~\ref{2-0}, we introduce the class of {\sl geometry-compatible conservation laws} on $\MM$
based on divergence-free flux, which we single out as a class of particular interest. We then define
the notions of entropy pair and entropy solution in the function space $L^\infty(\RR_+ \times \MM)$.
In Section~\ref{VA-0}, we show the existence of entropy solutions via the so-called vanishing diffusion method,
and establish that solutions to geometry-compatible conservation laws enjoy $L^p$ stability properties
(for all $p$) that are {\sl independent of the geometry} of the manifold.
Next, Section~\ref{3-0} is devoted to the general well-posedness theory on a Riemannian manifold,
and covers both $L^\infty$ solutions (following Kruzkov \cite{Kruzkov})
and measure-valued solutions (following DiPerna~\cite{DiPerna}).

In turn, we obtain a versatile framework
allowing us to establish the convergence toward the unique entropy solution
of any sequence of approximate solutions satisfying all of the entropy inequalities (possibly up to some small error).
Our geometry-independent bounds should be useful for the numerical analysis
of stable and robust, shock-capturing schemes adapted to hyperbolic equations posed on manifolds. In this direction,
the follow--up paper \cite{AmorimBenArtziLeFloch} will cover the convergence of finite volume
methods based on non-Cartesian meshes, and the derivation of total variation diminishing (TVD) estimates.

The remaining of the paper contains several generalizations of interest.
In Section~\ref{GE-0} we show that the well-posedness theory can be developed in the function space $L^1$
for {\sl general} conservation laws that need not be geometry-compatible; we prove that the $L^1$ contraction property
hold at this level of generality, but that for $p >1$ the $L^p$ stability properties may be violated.
Finally, in Section~\ref{LO-0} we discuss the well-posedness theory for conservation laws posed
on a Lorentzian manifold, including for instance the Schwarzschild spacetime.

\section{Preliminaries}
\label{PRE-0}

In this preliminary section, we briefly review some key concepts and results from Kruzkov's theory in the Euclidian space.
Consider the {\sl hyperbolic conservation law}
\be
\del_t u + \sum_{j=1}^n \del_j \big( f^j(u,\cdot) \big) = 0, \quad u=u(t,x) \in \RR, \, t>0, \, x \in \Rn,
\label{1.1}
\ee
where $\del_j := \del/\del x^j$  and
$f: \RR \times \Rn \to \Rn$ is a given smooth mapping, refered to as the {\sl flux} of \eqref{1.1}.
Observe that in \eqref{1.1} (as well as in \eqref{1.3} below) the divergence operator acts on the $x$-coordinate
and takes into account the dependence in $x=(x^j)$ of both the solution $u=u(t,x)$ and the flux
$f = (f^j(\ubar,x))$,  that is
$$
\aligned
& \sum_{j=1}^n \del_j \big( f^j(u(t,x),x) \big)
\\
& = \sum_{j=1}^n \Big( \big( \del_u f^j \big) (u(t,x),x) \, \del_j u(t,x) + \big(\del_j f^j\big)(u(t,x),x) \Big).
\endaligned
$$

As is well-known \cite{Dafermos,Hormander,Kruzkov,Lax,LeFloch,Volpert}, solutions of \eqref{1.1}
are generally discontinuous and must be understood in the weak sense of distributions. In addition,
for the sake of uniqueness one should restrict attention to {\sl entropy solutions} characterized by an infinite
family of inequalities. For simplicity, assume first that the flux is {\sl divergence-free,} that is
\be
\sum_{j=1}^n \big( \del_j f^j \big) (\ubar,x) = 0, \quad \ubar \in \RR, \, x \in \Rn.
\label{1.2}
\ee
(This assumption is motivated by Lemma~\ref{basicclaim} below.)
When the condition \eqref{1.2} is imposed, the {\sl entropy inequalities} associated with the equation \eqref{1.1}
take the form 
\be
\aligned
& \del_t U(u) + \sum_{j=1}^n \del_j \big( F^j(u,\cdot) \big) \leq 0,
\\
& U'' \geq 0,
\\
& F^j(\ubar,x) := \int^{\ubar} \del_u U(u') \, \del_u f^j(u',x) \, du', \quad x \in \Rn, \, \ubar \in \RR, \, j=1, \ldots, n.
\endaligned
\label{1.3}
\ee
Relying on Kruzkov's theory \cite{Kruzkov} and the arguments given later in this paper, one can check that,
given any divergence-free flux and initial data $u_0 \in L^1(\Rn) \cap L^\infty(\Rn)$,
the initial value problem \eqref{1.1},
\be
u(0,x) = u_0(x), \quad x \in \Rn,
\label{1.4}
\ee
admits a unique entropy solution $u \in L^\infty \bigl( \RR_+, L^1(\Rn) \cap L^\infty(\Rn) \bigr)$ which,
moreover, satisfies the {\sl $L^p$ stability properties}
\be
\|u(t)\|_{L^p(\Rn)} \leq \|u_0\|_{L^p(\Rn)}, \quad t \geq 0, \, p \in [1, \infty],
\label{1.5}
\ee
together with the {\sl $L^1$ contraction property:} for any two entropy solutions $u, v$
\be
\| v(t) - u(t) \|_{L^1(\Rn)} \leq \| v(t') - u(t')\|_{L^1(\Rn)}, \quad 0 \leq t' \leq t.
\label{1.6}
\ee
Furthermore, one has the following regularity property: If the initial data has bounded total variation, that is
$u_0 \in BV(\Rn)$,  then the solution has bounded variation for all times, that is
$u \in L^\infty \bigl( \RR_+, BV(\Rn))$.

In the Euclidian setting, the condition  \eqref{1.2} arises in the following way.

\begin{lemma}
\label{basicclaim}
A necessary and sufficient condition for
every {\sl smooth} solution of \eqref{1.1} to satisfy the infinite list of additional
{\rm conservation laws} \eqref{1.3} is that the divergence-free condition \eqref{1.2} holds.
\end{lemma}

\begin{proof}
This follows readily by multiplying the conservation law \eqref{1.1} by the derivative
$\del_u U(u)$ of an arbitrary function $U(u)$,
using the chain rule, and taking advantage of the fact that the function $U$ is arbitrary.
Namely, we have
$$
\aligned
\sum_{j=1}^n  \del_j (F^j(u,\cdot) )
& = \sum_{j=1}^n  (\del_u F^j)(u,\cdot) \, \del_j u + (\del_j F^j) (u,\cdot)
\\
& = \sum_{j=1}^n  \big( \del_u U \, \del_u f^j \big)(u,\cdot) \, \del_j u + (\del_j F^j) (u,\cdot),
\endaligned
$$
so that
$$
\del_t U(u) + \del_j \big( F^j(u,\cdot) \big)
=
\sum_{j=1}^n   (\del_j F^j)(u,\cdot) - \del_u U(u) \, (\del_j f^j) (u,\cdot).
$$
Now, by imposing that the right-hand side of the above identity vanishes for all entropy $U$
and by differentiating in $u$ the corresponding relation while using the definition of $F^j$ in \eqref{1.3},
we obtain
$$
\del_{uu} U(u)  \, (\del_j f^j) (u,\cdot) =0.
$$
Since $U$ is arbitrary, the desired claim follows.
\end{proof}

Observe that Kruzkov's theory in $\Rn$ also covers flux that need not satisfy our condition
\eqref{1.2}, and actually applies to {\sl balance laws} of the general form
\be
\del_t u + \sum_{j=1}^n \del_j \big( f^j(u,\cdot) \big) = h(u, \cdot), \quad u=u(t,x) \in \RR,
\label{1.balancelaw}
\ee
where $f=(f^j(\ubar,t,x))$ and $h=h(\ubar,t,x)$ are given smooth mappings. Provided $h$ is Lipschitz continuous in
the $u$-variable (uniformly in $t,x$), the initial value problem associated with \eqref{1.balancelaw}
admits a globally defined, unique, entropy solution
$u \in L^\infty_{\loc}( \RR_+ \times \Rn)$.  Note, however, that the general equation \eqref{1.balancelaw}
does not enjoy the $L^p$ stability, $L^1$ contraction and TVD properties.
Most of the literature on scalar conservation laws is restricted to the case of a
constant flux $f=f(\ubar)$ which, clearly, arises as a special case of \eqref{1.2}.
As pointed out earlier, there is no such concept as
a ``constant flux'' in the context of manifolds,
 and a suitable generalization of \eqref{1.2} will arise in our analysis.


\section{Geometry-compatible conservation laws}
\label{2-0}

\subsection*{Background and notation}

Throughout this paper, $(\MM,g)$ is a compact, oriented, $n$-dimensional Riemannian manifold.
As usual, the tangent space at a point $x \in \MM$ is denoted by
$T_x \MM$ and the tangent bundle by $T\MM : = \bigcup_{x \in \MM} T_x \MM$,
while the cotangent bundle is denoted by
$T^\star \MM = T_x^\star \MM$. The metric structure is determined by a positive-definite, $2$-covariant tensor field $g$, that
is, at each $x \in \MM$, $g_x$ is a inner product on $T_x\MM$.

In local coordinates $x=(x^j)$,
the derivations $\del_j := {\del \over \del x^j}$ form a basis of the tangent space $T_x \MM$, while
the differential forms $dx^j$ determine a basis of the cotangent space $T_x^\star \MM$.
Here and below, we use Einstein's summation convention on repeated indices so, for instance, in local coordinates
$$
g = g_{ij} \, dx^idx^j,  \qquad  g_{ij} = g(\del_i, \del_j).
$$
The notation $<X_x,Y_x>_{g_x} := g_x(X_x,Y_x)$ may also be used for the inner product of $X_x,Y_x \in T_x\MM$.

We denote by $(g^{ij})$ the inverse of the positive definite matrix $(g_{ij})$.
The metric tensor will be used to raise and lower indices, so that to each vector
$X=(X^j)$ one associates the covector $(X_j)$ via
$$
X_j := g_{ij} \, X^i.
$$

Recall that the differential $df$ of a function $f : \MM \to \RR$ is the section of the cotangent bundle
$T_x^\star \MM$, that is $df_x \in T_x^\star \MM$, defined by
$$
df_x (X_x) = X_x(f), \quad X_x \in T_x\MM.
$$
We denote by $dV_g$  the volume form on $\MM$, which in local coordinates reads
$$
dV_g = \sqrt{|g|} \, dx^1 \ldots dx^n,
$$
with $|g| := \det(g_{ij}) > 0$. The gradient of a function $h$ is obtained from
$dh$ by using the isomorphism between $T\MM$ and $T^*\MM$ via the Riemannian scalar
product.
In local coordinates, the gradient of a function $h$ is the vector field
$$
\grade h = (\nabla^j h) = {\del \over \del x^j} := g^{ij} \del_i h \, {\del \over \del x^j}.
$$

To the metric tensor $g$ one associates a  covariant derivative operator $\nabla$,
characterized by the condition $\nabla g = 0$. In particular, the covariant derivative
of a vector field $X = X^j \, {\del \over \del x^j}$ is the $(1,1)$-tensor field
whose coefficients in local coordinates are
$$
\nabla_j X^i := \del_j X^i + \Gamma_{kj}^i \, X^k, \quad i,j = 1, \ldots, n,
$$
where the Christoffel's symbols are determined from the metric tensor by
$$
\Gamma_{kj}^i := {1 \over 2} \, g^{il} \, \bigl( \del_k g_{lj} + \del_j g_{kj} - \del_l g_{kj} \bigr).
$$
The divergence of the vector field $X$ is the function
\be
\label{div}
\dive X := \nabla_j X^j = \del_j X^j + \Gamma_{kj}^j \, X^k.
\ee

Recall the duality relation
\be
\label{dual}
\int_\MM g \big( \grade h, X \big) \, dV_g = - \int_{\MM} h \, \dive X \, dV_g,
\ee
which is valid for any smooth function $h$ and vector field $X$, at least.
Interestingly enough, it follows from \eqref{dual} that, in local coordinates,
\be
\label{divloc}
\dive X = (\sqrt{|g|})^{-1} \, \del_j(\sqrt{|g|} \, X^j),
\ee
which shows that the divergence operator {\sl on vector fields} depends upon $|g|$, only.


\subsection*{A class of conservation laws}

The observation made earlier in Lem\-ma~\ref{basicclaim} in the Euclidian setting
provides us the motivation for the definition that we now introduce.
As was already pointed out in the introduction, although Kruzkov's theory applies to general
balance laws \eqref{1.balancelaw}, the subclass defined by \eqref{1.1} and the condition \eqref{1.2}
leads to many properties (maximum principle, $L^1$ contraction, $L^p$ stability) which
are very desirable features that we attempt to guarantee on the manifold $\MM$.
At this stage of the discussion, we restrict attention to smooth functions satisfying in a classical sense
the partial differential equations under consideration.

\begin{definition}
\label{div-cond}
1. A {\rm flux} on the manifold $\MM$ is a vector field $f=f_x(\ubar)$ depending upon the parameter $\ubar$
(the dependence in both variables being smooth).

2. The {\rm conservation law associated with the flux} $f$ on $\MM$ is \be \del_t u +
\dive (f(u))  = 0, \label{2.3} \ee where the unknown is the scalar-valued function
$u=u(t,x)$ defined for $t \geq 0$ and $x \in \MM$ and the divergence operator is
applied, for each fixed time $t$, to the vector field $x \mapsto f_x(u(t,x)) \in
T_x\MM$.

3. A flux is called {\rm geometry-compatible} if it satisfies the condition
\be
\dive f_x(\ubar) = 0, \qquad \ubar \in \RR, \, x \in \MM.
\label{2.divcondition}
\ee
\end{definition}

We will also refer to \eqref{2.3}-\eqref{2.divcondition} as a geometry-compatible conservation law.
Let us emphasize that \eqref{2.3} is a geometric partial differential equation which depends on the geometry
of the manifold, only, and is independent of any particular system of local coordinates on $M$.
In particular, all estimates derived for solutions of \eqref{2.3} must take a coordinate-independent form,
while, in the proofs, it will often be convenient to introduce a particular chart.

Observe that, by the divergence theorem on Riemannian manifolds
(\cite{Spivak}, for instance), all sufficiently smooth solutions of  \eqref{2.3} satisfy the
balance law
$$
{d \over dt} \int_S u(t,x) \, dV_g(x) = \int_{\del S} g_x \big( f_x(u(t,x)), \nu_x \big) \, dV_{\del S}(x),
$$
for every smooth, $n$-dimensional sub-manifold $S \subseteq \MM$ with boundary $\del S$.
Here, $\nu$ denotes the outward unit normal to $S$, and $V_{\del S}$ is the volume form induced on
the  boundary by restricting the metric $g$ to the tangent space $T(\del S)$.

In local coordinates, the flux components $f^j_x(\ubar)$ depend upon $\ubar,x$ and we will denote by
$\del_u f^j_x(\ubar)$ and $\del_k f^j_x(\ubar)$ its $\ubar$- and $x^k$-derivatives, respectively.
Before we proceed, observe that:

\begin{lemma} Let $f=f_x(\ubar)$ be a geometry-compatible flux on $\MM$.
In local coordinates the conservation law \eqref{2.3} takes the (nonconservative) form
\be
\del_t u(t,x) +  \big(\del_u f_x^j\big)(u(t,x)) \, \del_j u(t,x) = 0.
\label{2.6} \ee
\end{lemma}

\begin{proof} It follows from the definition of the divergence operator \eqref{div} that the conservation law \eqref{2.3}
can be written as
$$
\del_t u(t,x) + \del_j \big( f_x^j(u(t,x)) \big) + \Gamma_{kj}^j(x) \, f^k_x(u(t,x)) = 0.
$$
The function $\del_j \big( f^j_x(u(t,x) \big)$ is the sum of partial derivatives with respect to the local coordinate $x$
of the composite map $x \mapsto f_x(u(t,x))$, that is
$$
\del_j \big( f_x^j(u(t,x)) \big) =  \big(\del_u f_x^j\big)(u(t,x)) \, \del_j u(t,x) +  (\del_j f_x^j)(u(t,x)).
$$
On the other hand, the condition in Definition~\ref{div-cond} yields
$$
\big(\del_j f^j_x\big) (\ubar,x) + \Gamma_{kj}^j(x) \, f^k_x(\ubar) = 0, \quad \ubar \in \RR, \, x \in \MM.
$$
Writing this identity with $\ubar = u(t,x)$ and combining it with the above observations
lead us to \eqref{2.6}.
\end{proof}


\subsection*{Weak solutions}

We will be interested in measurable and bounded functions ${u \in L^\infty(\RR_+\times\MM)}$
satisfying \eqref{2.3} in the sense of distributions and
assuming a prescribed initial condition $u_0 \in L^\infty(\MM)$ at the time $t=0$:
\be
\label{2.init} u(0,x) = u_0(x), \quad x \in \MM.
\ee
For the sake of uniqueness, we
need an analogue of the entropy inequalities \eqref{1.3} proposed by Kruzkov in the
Euclidian space $\RR^n$. These inequalities are derived as follows.

First of all, consider smooth solutions $u=u(t,x)$ of \eqref{2.3} and multiply
the conservation law by the derivative $\del_u U(u)$ of an arbitrary function $U$. We obtain
$$
\del_t U(u(t,x)) + \del_u U(u(t,x)) \, \dive \bigl( f_x(u(t,x)) \bigr) = 0,
$$
which suggests to search for a vector field $F=F_x(\ubar)$ such that
$$
\del_u U(u(t,x)) \, \dive \big( f_x(u(t,x)) \big) = \dive \big( F_x(u(t,x)) \big).
$$
This relation should hold for all functions $u=u(t,x)$, and is equivalent to impose that the flux
components $f^j_x(\ubar)$ satisfy the two partial differential equations
\be
\label {fF}
\aligned
& \del_u U(\ubar) \, \del_u f^j_x(\ubar) =  \del_u F^i_x(\ubar), \qquad 1 \leq j \leq n,
\\
& \del_u U(\ubar) \, \big(\del_j f^j_x(\ubar) + \Gamma_{kj}^j(x) \, f^k_x(\ubar) \big)
=  \del_j F^j_x(\ubar) + \Gamma_{kj}^j(x) \, F^k_x(\ubar),
\endaligned
\ee for all $\ubar \in \RR$ and $x \in \MM$. The first equation implies that $F_x$ is
given by \be \label{ent-flux1} F_x(\ubar) := \int^{\ubar} \del_w U(w) \, \del_w f_x(w)
\, dw \in T_x\MM. \ee Differentiating the second equation with respect to $\ubar$ and
using \eqref{ent-flux1} we find
$$
\del_{uu} U(\ubar) \, \big(\del_j f^j_x(\ubar) + \Gamma_{kj}^j(x) \, f^k_x(\ubar) \big) = 0.
$$
Since $U$ is arbitrary, a necessary and sufficient condition for \eqref{fF} to
hold is that $f_x(\ubar)$ satisfies \eqref{2.divcondition}.
Note in passing that necessarily $F_x(\ubar)$ satisfies also this condition.

This derivation provides us with an analogue of the definition \eqref{1.3} of the Euclidian case.
In turn, smooth solutions of any geometry-compatible conservation law \eqref{2.3}
satisfy the additional conservation laws
\be
\del_t U(u(t,x)) + \dive \big( F_x(u(t,x)) \big) = 0,
\label{2.entropyequality}
\ee
where $U$ is arbitrary and $F$ is given by \eqref{ent-flux1}.

We are now in position to introduce the notion of entropy solution which generalizes Kruzkov's notion
to the case of manifolds. We require the equalities \eqref{2.entropyequality} to hold
as inequalities only, as this arises naturally by the vanishing diffusion method. (See~Section~\ref{VA-0}).

\begin{definition}
\label{2-3} Let $f=f_x(\ubar)$ be a geometry-compatible flux on the Riemannian manifold $(M,g)$.

1. A {\rm convex entropy/entropy-flux pair} is a pair $(U,F)$ where $U: \RR \to \RR$
is a convex function, and $F=F_x(\ubar)$ is the vector field depending on the
parameter $\ubar$ defined by \be \label{ent-flux} F_x(\ubar) := \int^{\ubar}
\hskip-.15cm \del_{u'} U(u') \, \del_{u'} f_x(u') \, du',  \qquad \ubar \in \RR, \, x
\in \MM. \ee

2. Given a function $u_0 \in L^\infty(\MM)$, a function $u$ in $L^\infty\bigl(\RR_+, L^\infty(\MM)\bigr)$ is
called an {\rm entropy solution} to the initial value problem \eqref{2.3}-\eqref{2.init} if
the following {\rm entropy inequalities} hold
\be
\aligned
& \iint_{\RR_+ \times \MM}  \hskip-.15cm \Bigl( U(u(t,x)) \, \del_t \theta(t,x)
+ g_x\big(F_x(u(t,x)), \grade \theta(t,x) \big) \Bigr) \, dV_g(x) dt
\\
& + \int_\MM  \hskip-.15cm U(u_0(x)) \, \theta(0,x) \, dV_g(x) \geq 0,
\endaligned
\label{ent-sol}
\ee
for every convex entropy/entropy flux pair $(U,F)$
and all smooth functions $\theta = \theta(t,x) \geq 0$ compactly supported in $[0,\infty)\times \MM$.
\end{definition}


\section{Vanishing diffusion method on manifolds}
\label{VA-0}

In the present section, we construct solutions of the conservation law \eqref{2.3} via the vanishing diffusion
method. We establish the existence of smooth solutions to a regularized version of \eqref{2.3} taking into account
a small diffusion term, and then establish that these solutions converge to an entropy solution
in the sense of Definition~\ref{2-3} as the diffusion tends to zero.
For convenience, we provide a proof that is based on a uniform total variation bound,
and refer the reader to the forthcoming section for an alternative
approach based on the concept of measure-valued solution.

Recall that $\MM$ endowed with the metric tensor $g$
is a compact, oriented, $n$-dimensional Riemannian manifold.
We denote by $L^p(\MM; dV_g)$ the usual Lebesgue spaces on the manifold $(\MM,g)$.
The Sobolev space $H^1(\MM; dV_g)$ is the space of all functions $h \in L^2(\MM; dV_g)$
such that $\grade h \in L^2(\MM; dV_g)$ ---where the gradient is defined in the distributional sense,
via the formula \eqref{dual}. The spaces $H^k(\MM; dV_g)$ for $k\geq 2$ are defined similarly.

The {\sl total variation} of a  function $h \in L^1(\MM; dV_g)$ is defined by
$$
TV(h) : = \sup_{\|\psi\|_{L^\infty(\MM)} \leq 1} \int_\MM h \, \dive \psi \, dV_g,
$$
where $\psi$ describes all $C^1$ vector fields on the manifold. We denote by $BV(\MM;
dV_g) \subset L^1(\MM; dV_g)$ the space of all functions $h$ with finite total variation,
endowed with the norm $\|h\|_{L^1(\MM; dV_g)} +TV(h)$. For material on BV
functions we refer to \cite{Federer,Volpert}. In particular, it is well-known that the
inclusion map of $BV(\MM; dV_g)$ in $L^1(\MM; dV_g)$ is compact.

It is clear that, for all smooth functions $h : \MM \to \RR$,
$$
TV(h) = \int_\MM |\grade h|_g \, dV_g,
$$
where $| \cdot |_g$ denotes the Riemannian norm associated with $g$.
Observe also that, by taking a partition of unity determining a finite covering of $\MM$ by coordinate patches
we can easily ``localize'' the concept of total variation on the manifold, as follows.

\begin{lemma}
\label{VA.XTV}
There exist finitely many (smooth) vector fields $X^{(1)}$,\ldots, $X^{(L)}$ on $\MM$,
a constant $C_0 \geq 1$ (depending only on the manifold), and a chart covering the manifold such that:
\par
(ii)  each vector field $X^{(l)}$ is supported in a coordinate patch and
$$
\text{Span} \big\{X_x^{(1)},...,X_x^{(L)} \big\} : = T_x\MM, \qquad x \in \MM;
$$
(ii) and for every smooth function $h:\MM \to \RR$
$$
{1 \over C_0} \, \sum_{l=1}^L \int_\MM |X^{(l)} (h)|_g \, dV_g
\leq
TV(h) \leq C_0 \, \sum_{l=1}^L \int_\MM |X^{(l)} (h)|_g \, dV_g.
$$
\end{lemma}


\

An initial data $u_0$ being given in $L^\infty(\MM) \cap BV(\MM; dV_g)$, we want to
find a solution $u^\eps=u^\eps(t,x)$ to the initial value problem \be \del_t u^\eps  +
\dive \big( f_x(u^\eps) \big) = \eps \, \Delta_g u^\eps, \qquad x \in \MM, \, t \geq
0, \label{VA.1} \ee \be u^\eps(0,x) = u_0^\eps(x), \quad x \in \MM, \label{VA.2} \ee
where $\Delta_g$ denotes the Laplace operator on the manifold $\MM$,
$$
\aligned \Delta_g v & := \dive \grade v
\\
& = g^{ij}\, \big( {\del^2 v \over \del x^i \del x^j} - \Gamma^k_{ij} \, {\del v \over \del x^k } \big).
\endaligned
$$
In \eqref{VA.2}, $u_0^\eps: \MM \to \RR$ is a sequence of smooth functions satisfying
\be
\aligned
& \|u_0^\eps\|_{L^p(\MM)} \leq \|u_0\|_{L^p(\MM)}, \qquad p \in [1,\infty],
\\
& TV(u_0^\eps) \leq TV(u_0),
\\
 & \eps \, \|u_0^\eps\|_{H^2(\MM; dV_g)} \leq C\,TV(u_0), \quad \text {for some $C>0$
 depending only on $\MM$,} 
\\
& u_0^\eps \to u_0 \quad \text{ a.e. on }  \MM.
\endaligned
\label{VA.initreg}
\ee

In the following, we use the notation $u(t)$ instead of $u(t,\cdot)$. We begin with:

\begin{theorem} {\rm (Regularized conservation law.)}
\label{VA.viscous-M} Let ${f=f_x(\ubar)}$ be a geometry-compatible flux on a
Riemannian manifold $(\MM,g)$. Given $u_0^\eps \in C^\infty(\MM)$ 
satisfying \eqref{VA.initreg} there exists a unique solution $u^\eps \in
{C^\infty(\RR_+ \times \MM)}$ to the initial value problem \eqref{VA.1}-\eqref{VA.2}.
Moreover, for each ${1 \leq p \leq \infty}$ the solution satisfies \be
\|u^\eps(t)\|_{L^p(\MM; dV_g)} \leq \|u^\eps(t')\|_{L^p(\MM; dV_g)}, \qquad 0 \leq t'
\leq t \label{VA.3} \ee and, for any two solutions $u^\eps$ and $v^\eps$, \be
\|v^\eps(t) - u^\eps(t) \|_{L^1(\MM; dV_g)} \leq \| v^\eps(t') - u^\eps(t')
\|_{L^1(\MM;dV_g)}, \quad 0 \leq t' \leq t. \label{VA.4} \ee In addition, for every
convex entropy/entropy flux pair $(U,F_x)$ the solution $u^\eps$ satisfies the entropy
inequality \be \del_t U(u^\eps)  + \dive \big( F(u^\eps) \big) \leq \eps \, \Delta_g
U(u^\eps). \label{VA.5} \ee
\end{theorem}

All proofs are postponed to the end of this section.

Our next goal is to prove the strong convergence of the family $\{u^\eps\}_{\eps>0}$ to an entropy solution of
\eqref{2.3}-\eqref{2.init}, as $\eps\rightarrow 0$. We will assume here that the initial data has bounded total variation,
and we refer to the forthcoming section for a more general argument which does not require this assumption.
From the estimates for the solutions obtained in Theorem~\ref{VA.viscous-M} we can deduce a uniform
total variation bound in space and time.

\begin{lemma} {\rm (BV bounds for diffusive approximations.)}
\label{VA.viscous-VAR}
There exists a positive constant $C_1$ depending on the geometry of $\MM$
and $\| u_0\|_{L^\infty(\MM)}$, only, such that the solutions $u^\eps$
given by Theorem~\ref{VA.viscous-M} satisfy
\be
\label{VA.TV-estimate}
TV(u^\eps(t)) \leq e^{C_1 \, t} \, \big( 1 + TV(u_0) \big), \qquad t \in \RR_+,
\ee
\be
\label{VA.t-cont}
\| \del_t u^\eps(t)\|_{L^1(\MM; dV_g)}
\leq  TV(u_0) + \eps \, \|\Delta_g u_0^\eps\|_{L^2(\MM; dV_g)},
 \quad t \in \RR_+.
\ee
\end{lemma}

Observe that the right-hand side of \eqref{VA.TV-estimate} is independent of $\eps$,
while the term $\eps \, C_1 \|\Delta_g u_0^\eps\|_{L^2(\MM)}$ arising in
\eqref{VA.t-cont} is dominated by $C\,TV(u_0)$ thanks to our assumption  \eqref{VA.initreg}.
At this juncture, we refer to 
\cite{AmorimBenArtziLeFloch} in which a sharp version of the above estimates is derived with a
specific constant $C_1$ depending on the geometry of the manifold.

We are now in a position to state and prove the main result of this section.

\begin{theorem} {(\rm Convergence of the vanishing diffusion method on a Riemannian manifold.)}
\label{VA.entropy-sol}
Let $f=f_x(\ubar)$ be a geometry-compatible flux on a Riemannian manifold $(\MM,g)$.
Given any initial condition $u_0\in L^\infty(\MM) \cap BV(\MM; dV_g)$ there exists an entropy solution
$u \in L^\infty(\RR_+ \times \MM)$ to the initial value problem \eqref{2.3}-\eqref{2.init}
in the sense of Definition~\ref{2-3}, which is the limit of the sequence $u^\eps$ constructed in
Theorem~\ref{VA.viscous-M} by vanishing diffusion.

Moreover, this solution has the following properties:
$$
\|u(t) \|_{L^p(\MM; dV_g)} \leq \|u_0\|_{L^p(\MM;dV_g)}, \qquad t \in \RR_+, \, 1 \leq p \leq \infty,
$$
and, for some constant $C_1>0$ depending on  $\| u_0\|_{L^\infty(\MM)}$ and the geometry of $\MM$ only,
\be
\aligned
& TV(u(t)) \leq e^{C_1 \, t} \, ( 1 + TV(u_0) ), \qquad t \in \RR_+,
\\
& \| u(t) - u(t') \|_{L^1(\MM;dV_g)} \leq TV(u_0) \, |t - t'|, \quad  0 \leq t' \leq t.
\endaligned
\label{VA.boundsonu} \ee Furthermore, if $u,v$ are entropy solutions associated with
some initial data $u_0, v_0$, respectively, it hold \be \label{contract} \| v(t) -
u(t) \|_{L^1(\MM; dV_g)} \leq \| v_0 - u_0 \|_{L^1(\MM; dV_g)}, \qquad t \in \RR_+.
\ee
\end{theorem}


\begin{proof}[Proof of Theorem~\ref{VA.viscous-M}]
By introducing a local chart of coordinates and using the condition \eqref{2.divcondition},
one can reduce the equation \eqref{VA.1} to a parabolic equation in the Euclidian space $\Rn$:
\be
\del_t u^\eps + (\del_u  f^j_x)(u^\eps) \, \del_j u^\eps
= \eps \, g^{ij} \, \big( \del_i \del_j u^\eps - \Gamma^k_{ij} \, \del_k u^\eps \big).
\label{VA.visc-eq}
\ee
The {\sl local in time} existence, uniqueness, and regularity of the
solution $u^\eps$ follows as in the Euclidian setting. By putting together several charts to cover the whole manifold,
one arrives at the {\sl local in time} well-posedness for the diffusive conservation law \eqref{VA.1}.
Next, in order to extend the solution to arbitrary times, one needs a uniform estimate on the sup norm of the solution.
This indeed follows in the form stated in \eqref{VA.3} for $p=\infty$, from the standard maximum principle for parabolic equations.

Next, in every strip $[0,T] \times \MM$ the dual equation associated with \eqref{VA.1} reads
$$
\del_t \varphi + g_x \big( f_x(u^\eps), \grade \varphi\big) + \eps \, \Delta_g \varphi = 0,
$$
for which a ``terminal value" problem with data given at $t=T$ is now considered.
This equation satisfies the maximum principle and, by duality, the inequality \eqref{VA.3} for $p=1$
follows. The inequality for intermediate values $p\in (1,\infty)$ is then obtained by straightforward interpolation.

We next proceed with the derivation of \eqref{VA.4}. The function $w^\eps:=v^\eps - u^\eps$
satisfies the equation
\be
\label{VA.weps}
\del_t w^\eps + \dive (b^\eps_x \, w^\eps) = \eps \, \Delta_g w^\eps,
\ee
where
$$
b_x^\eps : = {f_x(v^\eps) - f_x(u^\eps)  \over v^\eps - u^\eps}
$$
is a smooth vector field on $\MM$. As in the argument above, the dual equation associated with \eqref{VA.weps}
satisfies the maximum principle and so, by duality,
$$
\|w^\eps(t)  \|_{L^1(\MM; dV_g)} \leq \| w^\eps(t')  \|_{L^1(\MM;dV_g)}, \quad 0 \leq t' \leq t,
$$
which is precisely \eqref{VA.4}.

Finally, to show \eqref{VA.5} we multiply \eqref{VA.1} by $\del_u U(u^\eps)$, where $(U, F)$
is a convex entropy pair in the sense of Definition~\ref{2-3}. We get
$$
\del_t U(u^\eps) + \dive F_x(u^\eps) = \eps \, \del_u U(u^\eps) \, \Delta_g u^\eps.
$$
In view of the definition of $\Delta_g$ we can write in local coordinates
\be
\label{equa7}
\del_u U(u^\eps) \, \Delta_g u^\eps = \Delta_g U(u^\eps) - \del_u^2 U(u^\eps) \, g^{jk} \, \del_j u^\eps \del_k u^\eps.
\ee
Since $\del_u^2 U$ is  non-negative ($U$ being convex) and the matrix $(g^{jk})$ is
positive definite (since $g$ is a Riemannian metric tensor), the second term in the right-hand side of
\eqref{equa7} is non-positive. This establishes the desired inequality \eqref{VA.5}.
\end{proof}

\begin{proof}[Proof of Lemma~\ref{VA.viscous-VAR}]
To establish \eqref{VA.TV-estimate} we rely on Lemma~\ref{VA.XTV}: 
without loss of generality, it is sufficient to estimate $\|\del_1(\psi
u^\eps(t))\|_{L^1(\MM; dV_g)}$ for every smooth $\psi$ supported in a
coordinate patch $\Ocal$ with $\|\psi\|_{L^\infty(\MM)} \leq 1$. 
Multiplying  \eqref{VA.visc-eq} by $\psi$ and differentiating with respect to $x^1$, the function
$w^\eps := \del_1(\psi \,u^\eps)$ is found to satisfy
$$
\aligned
& \del_t w^\eps(t,x) + \del_j \big( (\del_u f^j_x)(u^\eps(t,x)) \, w^\eps(t,x)\big)
\\
& = A(x,u^\eps, \grade u^\eps)(t,x) + \eps \, \Delta_g w^\eps(t,x),
\endaligned
$$
where $A$ depends linearly on $\grade u^\eps$ and smoothly upon $x, u^\eps$ (and is supported in $\Ocal$).
We now multiply the identity above by $\sgn(w^\eps)$ and integrate over $\MM$.
At this juncture we  note that, since $w^\eps\in H^2(\MM; dV_g)$ at least,
$$
\aligned
&  \sgn (w^\eps) \, \del_t w^\eps = \del_t |w^\eps|,
\qquad \sgn (w^\eps) \, \del_j w^\eps = \del_j |w^\eps|,
\\
& \sgn (w^\eps) \, \Delta_g w^\eps \leq \Delta_g |w^\eps|,
\endaligned
$$
where the latter inequality must be understood in the sense of distributions.

Relying on the bound $\| u^\eps(t)\|_{L^\infty(\MM)} \leq \| u_0\|_{L^\infty(\MM)}$ which follows
from \eqref{VA.3}, we then deduce that there exists a constant $C_1>0$ (depending on both
$\| u_0\|_{L^\infty(\MM)}$ and the geometry of $\MM$) such that
$$
{d \over dt} \int_\MM |w^\eps(t)| \, dV_g
\leq
C_1 + C_1 \, \int_\MM |\grade u^\eps(t)|_g \, dV_g.
$$
This estimate can now be repeated for each of the (finitely many) vector fields $X^{(1)},\ldots,X^{(L)}$, and
the desired conclusion \eqref{VA.TV-estimate} finally follows from Gronwall's lemma.

To establish \eqref{VA.t-cont} we differentiate \eqref{VA.visc-eq} with respect to $t$ and note that
$$
\aligned
& \dive\big( \partial_tf_x(u^\eps)\big) = \partial_t\dive\big(f_x(u^\eps)\big),
\\
& \Delta_g \big(\partial_tf_x(u^\eps) \big) = \partial_t\Delta_g \big( f_x(u^\eps)\big).
\endaligned
$$
Also, the vector field $\partial_tf_x(u^\eps(t,x))$ satisfies the condition
\eqref{2.divcondition} (with respect to the explicit dependence on $x$).

It therefore follows that the function $z^\eps = \del_t( u^\eps)$ satisfies \be
\label{zeps}
 \del_t z^\eps(t,x) + \frac{\del}{\del y_j} \big(  (\del_u f^j_x)(u^\eps(t,y)) \, z^\eps(t,y)
\big)_{y=x} = \eps \, \Delta_g z^\eps (t,x).
\ee
(In other words, the explicit dependence of $(\del_u f^j_x)(u^\eps(t,x))$ on $x$ is not differentiated.)
Note that
$$
\aligned 
& \sgn(z^\eps)\frac{\del}{\del y_j} \big(  (\del_u f^j_x)(u^\eps(t,y)) \, z^\eps(t,y)
\big)_{y=x}
\\
& =\frac{\del}{\del y_j} \big(  (\del_u f^j_x)(u^\eps(t,x)) \, |z^\eps(t,x)|
\big)_{y=x}.
\endaligned 
$$
Defining the vector field $r_x:=(\del_u f_x)(u^\eps(t,x)) \,
|z^\eps(t,x)|$ we have, by the condition \eqref{2.divcondition}, that
the second term in the left-hand side of \eqref{zeps} satisfies
$$
\del_j \big(  (\del_u f^j_x)(u^\eps(t,x)) \,
|z^\eps(t,x)| \big)=\dive r_x.
$$

Multiplying \eqref{zeps} by $\sgn(z^\eps)$ and integrating over $\MM$, 
we obtain
$$
{d \over dt} \int_\MM |z^\eps(t)| \, dV_g \leq 0. 
$$
In turn, we can write
$$
\int_\MM |z^\eps(t)| \, dV_g \leq \int_\MM |z^\eps(0)| \, dV_g.
$$
Alternatively, this inequality can be derived from \eqref{VA.4} by choosing $v(t,x) = v(t+\alpha,x)$ with $\alpha \to 0$.
Using the equation \eqref{VA.1} we can estimate the above term at $t=0$:
$$
\aligned
\int_\MM |z^\eps(0)| \, dV_g
& \leq \|\grade u_0^\eps\|_{L^1(\MM; dV_g)} + \eps \, \| \Delta_g u_0^\eps\|_{L^2(\MM; dV_g)}
\\
& \leq TV( u_0^\eps) + \eps \, \| \Delta_g u_0^\eps\|_{L^2(\MM; dV_g)}.
\endaligned
$$
Finally, applying this estimate to each of the vector fields $X^{(1)},\ldots,X^{(L)}$,
we arrive at the estimate \eqref{VA.t-cont}. This completes the proof of Lemma~\ref{VA.viscous-VAR}.
\end{proof}

\begin{proof}[Proof of Theorem~\ref{VA.entropy-sol}]
In view of the sup norm estimate \eqref{VA.3} and the uniform, space $L^1$ estimates
\eqref{VA.TV-estimate}, we see that the sequence $u^\eps(t)$ is uniformly bounded in
$BV(\MM; dV_g)$ for every time $t$ and is therefore compact in $L^1(\MM; dV_g)$.
Applying this argument at all rational times $t$ and then picking up a diagonal
sequence we can ensure that a subsequence $u^{\eps_j}(t)$ converges to some limit
$u(t)$ in the $L^1$ norm, as $\eps_j\rightarrow 0,$ for all rational $t.$ Next, in
view of the uniform time estimate \eqref{VA.t-cont}, the limit $u(t)$ extends to all
values of the time variable, with $u^{\eps_j}(t) \to u(t)$ in $L^1(\MM; dV_g)$.
Letting $\eps_j \to 0$ in the inequalities \eqref{VA.5} one then deduces that $u$
satisfies all of the entropy inequalities
$$
\del_t U(u) + \dive \big( F(u) \big) \leq 0
$$
in the weak sense \eqref{ent-sol}.

In fact, by working in  a localized coordinate
patch and using the entropy formulation ~\eqref {VA.5} with $U(u)=|u-k|,$
$F_x(u)=\sgn(u-k)(f_x(u)-f_x(k)),$ one can repeat the Kruzkov classical theorem in
order to obtain a "localized" $L^1$ contraction property. By patching together
(finitely many) coordinate patches we obtain a global estimate of the form \be  \|
v(t) - u(t) \|_{L^1(\MM; dV_g)} \leq C\| v_0 - u_0 \|_{L^1(\MM; dV_g)}, \qquad t \in
[0,T], \ee
  where $C>0$ depends on $T$ (and $\MM).$ This is sufficient to imply uniqueness of
  the entropy solution $u(t)$, hence the convergence of the whole family $u^\eps.$
  At this point we can invoke the parabolic estimate ~\eqref{VA.4} and conclude that
  $C=1,$ thus establishing the $L^1$ contraction property \eqref{contract}.
 All the other  estimates now follow from Theorem~\ref{VA.viscous-M} and
Lemma~\ref{VA.viscous-VAR}.
\end{proof}


\section{The well-posedness theory}
\label{3-0}

Theorem~\ref{VA.entropy-sol} assumes that the initial data has bounded variation and
provides the existence of (locally BV) entropy solutions constructed by vanishing diffusion.
In the present section, we obtain a generalization to $L^\infty$ initial data and also establish
the uniqueness of the entropy solution in this larger class.
We emphasize that the class of $L^\infty$ solutions is completely
natural for geometry-compatible conservation laws, since the $L^\infty$
estimate is independent of the geometry while the BV estimate \eqref{VA.boundsonu}, in general, depends upon it.

Our generalization of Kruzkov's theory to manifolds is as follows.

\begin{theorem} {\rm (Well-posedness theory in $L^\infty$
for geometry-compatible conservation laws.)}
\label{3-1}
Let $f=f_x(\ubar)$ be a geometry-compatible flux on a compact, oriented, Riemannian manifold $(\MM,g)$.
Given any initial data $u_0 \in L^\infty(\MM)$ there exists a unique entropy solution
$u \in L^\infty(\RR_+ \times \MM)$ to the initial value problem \eqref{2.3}-\eqref{2.init}
in the sense of Definition~\ref{2-3}. Moreover, for each $1 \leq p \leq \infty$ the solution satisfies
\be
\|u(t) \|_{L^p(\MM; dV_g)} \leq \|u_0\|_{L^p(\MM;dV_g)}, \qquad t \in \RR_+,
\label{3.1}
\ee
and, given any two entropy solutions $u,v$ associated with some initial data $u_0, v_0$, respectively,
\be
\| v(t) - u(t) \|_{L^1(\MM; dV_g)} \leq \| v_0 - u_0 \|_{L^1(\MM; dV_g)}, \qquad t \in \RR_+.
\label{3.2}
\ee
\end{theorem}

Furthermore, we will see in the proof that the following inequality holds in the sense
of distributions \be \del_t | v - u | + \dive \big( \sgn(u - v) \, (f_x(v) - f_x(u)
)\big) \leq 0. \label{3.3} \ee

Following DiPerna~\cite{DiPerna} we introduce the (larger) class of entropy measure-valued solutions
and, in fact, establish a much stronger version of Theorem~\ref{3-1}. As observed in \cite{DiPerna},
Kruzkov's arguments take a simpler form in the measure-valued setting.
We consider {\sl measure-valued maps} $\nu = \nu_{t,x}$, that is,
weakly measurable mappings $(t,x) \in \RR_+ \times \MM \mapsto \nu_{t,x}$
taking their values in the space of probability measures on $\RR$
with support included in a fixed compact interval of $\RR$. The action of the measure $\nu$ on a function $U$ will be denoted
by
$$
\la \nu_{t,x}, U \ra := \int_\RR U(\ubar) \, d\nu_{t,x}(\ubar).
$$
The weak measurability property means that the map $\la \nu_{t,x}, U \ra$ is measurable in $(t,x)$ for each $U$.

\begin{definition}
\label{3-3} Let $f=f_x(\ubar)$ be a geometry-compatible flux on a Riemannian manifold
$(\MM,g)$. Given any initial condition $u_0 \in L^\infty(\MM)$, a measure-valued map
$(t,x) \in \MM\times \RR_+ \mapsto \nu_{t,x}$ is called an {\rm entropy measure-valued
solution} to the initial value problem  \eqref{2.3}-\eqref{2.init} if, for every
convex entropy/entropy flux pair $(U, F_x)$ (see \eqref{ent-flux}), \be \aligned &
\iint_{\RR_+ \times \MM}  \hskip-.15cm \Bigl( \la \nu_{t,x}, U \ra \, \del_t
\theta(t,x) + g_x\big(\la \nu_{t,x}, F_x \ra, \grade \theta(t,x) \big) \Bigr) \,
dV_g(x) dt
\\
& + \int_\MM  \hskip-.15cm U(u_0(x)) \, \theta(0,x) \, dV_g(x) \geq 0,
\endaligned
\label{3.5}
\ee
for every smooth function $\theta=\theta(t,x) \geq 0$ compactly supported in
$[0,+\infty) \times \MM$.
\end{definition}

We will now prove that:

\begin{theorem} {\rm (The well-posedness theory in the measure-valued class for geometry-compatible conservation laws.)}
\label{3-4}
Let $f=f_x(\ubar)$ be a geo\-metry-compatible flux on a compact, oriented, Riemannian manifold $(\MM,g)$.
Let $u_0$ be in $L^\infty(\MM)$ and $\nu$ be an entropy measure-valued solution (in the sense of Definition~\ref{3-3})
to the initial value problem \eqref{2.3}-\eqref{2.init}. 
Then, for almost every $(t,x)$, the measure $ \nu_{t,x}$ is a Dirac mass,
i.e.~of the form
$$
\nu_{t,x} = \delta_{u(t,x)},
$$
where the function $u \in L^\infty(\RR_+ \times \MM)$ is the unique entropy solution
to the problem  \eqref{2.3}-\eqref{2.init} in the sense of Definition~\ref{2.3}. Moreover, the solution
satisfies the properties \eqref{3.1}--\eqref{3.3}, and the initial data is assumed in the strong sense
\be
\limsup_{t \to 0+} \int_\MM |u(t,x) - u_0(x)| \, dV_g(x) = 0.
\label{A-strong-initial}
\ee
\end{theorem}

We have already shown in Section~\ref{VA-0} that, provided the initial data have
bounded variation, entropy solutions can be constructed via vanishing diffusion. The
proof relies on the compactness the inclusion of $BV$ into $L^1$. Thanks to
Theorem~\ref{3-4}, we can now provide an alternative, more general proof which is
valid for $L^\infty$ initial data. In view of the uniform $L^\infty$ estimate (i.e.,
the maximum principle which only requires the initial data to be in $L^\infty$),
vanishing diffusion approximations $u^\eps$ generate a measure-valued solution $\nu$.
Moreover, $\nu$ can be easily checked to satisfy all of the entropy inequalities and,
in turn, Theorem~\ref{3-1} follows from Theorem~\ref{3-4}.

Finally, we also have the following generalization of Theorem~\ref{VA.entropy-sol}.

\begin{corollary} The convergence result of the vanishing diffusion approximations
in Theorem~\ref{VA.entropy-sol} remains valid if the initial data $u_0$ is solely in
$L^\infty(\MM)$. All estimates therein, except for the total variation estimate, still
hold.
\end{corollary}

In the course of proving Theorem~\ref{3-4}, we will need some a~priori regularity of measure-valued solutions,
especially the fact that the initial data $u_0$ is automatically assumed in a strong sense.

\begin{lemma}
\label{3-5} Let $\nu=\nu_{t,x}$ be an entropy measure-valued solution of
\eqref{2.3}-\eqref{2.init}, where it is assumed that $u_0 \in L^\infty(\MM)$ and the
supports of the measures $\nu_{t,x}$ are all included in a fixed common interval.
Then:

(a) For every convex function $U=U(u)$ and every smooth function $\theta=\theta(x)$, the function
\be
t \mapsto \int_\MM \la \nu_{t,x}, U \ra \, \theta \, dV_g(x)
\label{A-2.10}
\ee
has locally bounded total variation and admits a trace as $t \to 0+$.

(b)  For every function $U=U(u, x)$, which is convex
in $u$, measurable in $x$, and such that $|U(u,x)| \leq c \, |u| + |\widetilde U(x)|$
with $\widetilde U \in L^1(\MM; dV_g)$ and $c \geq 0$,
\be
\limsup_{t \to 0+} \int_\MM \la \nu_{t,x}, U(., x)\ra \, dV_g(x)
\leq
\int_\MM U(u_0(x), x) \, dV_g(x).
\label{A-2.11}
\ee

(c)  In particular, $\nu$ assumes its initial data $u_0$
in the following strong sense:
\be
\limsup_{t \to 0+} \int_\MM \int_\RR  |\ubar - u_0(x)| \, d\nu_{t,x}(\ubar) \, dV_g(x) = 0.
\label{A-2.12}
\ee
\end{lemma}

 
\begin{proof}[Proof of Lemma~\ref{3-5}]
Using in the weak formulation ~\eqref{3.5} a function $\theta(t,x) = \theta_1(x) \,
\theta_2(t)$, compactly supported in $[0, \infty) \times \MM$ and having $\theta_1,
\theta_2 \geq 0$, we obtain
$$
\aligned
& \int_0^\infty {d\theta_2 \over dt} \int_\MM \la \nu, U \ra \, \theta_1 \, dV_g(x) dt
  + \theta_2(0) \, \int_\MM U(u_0) \, \theta_1 \, dV_g(x)
\\
& \geq
  - \int_0^\infty \theta_2 \int_\MM g_x \big( \grade \theta_1, \la \nu, F_x \ra \bigr) \, dV_g(x) dt
\\
& \geq - C_1 \, \int_0^\infty \theta_2 \, dt,
\endaligned
$$
for some constant $C_1>0$ depending on $\theta_1$ (and the common support of
$\nu_{t,x}$). Thus the function
$$
V_1(t) := - C_1 \, t +  \int_\MM \la \nu_{t,x}, U \ra \, \theta_1 \, dV_g(x)
$$
satisfies the inequality
\be
- \int_0^\infty V_1(t) \, {d\theta_2\over dt} \, dt
\leq
\theta_2(0) \, \int_\MM U(u_0) \, \theta_1 \, dV_g(x).
\label{A-2.13}
\ee
Using a test-function $\theta_2 \geq 0$
compactly supported in $(0,\infty)$, we find
$$
- \int_0^\infty V_1(t) \, {d\theta_2\over dt} \, dt \leq 0.
$$
That is, the function $V_1(t)$ is decreasing
and, therefore, has locally bounded total variation. Since it is uniformly bounded,
$V_1(t)$ has a limit as $t \to 0+$. This proves (a).

To establish the item (b), we fix a time $t_0 > 0$ and consider the sequence of continuous functions
$$
\theta^\eps_2(t) = \begin{cases}
1,                  &  t \in [0, t_0],
\\
(t_0 + \eps -t)/\eps,   &  t \in [t_0, t_0+\eps],
\\
0,                  &  t \geq t_0+\eps.
\end{cases}
$$
Relying on the property (a) above, we see that
$$
- \int_0^\infty V_1(t) \, {d\theta_2^\eps\over dt} \, dt
\to
V_1(t_0+).
$$
Since $\theta_2^\eps(0) = 1$ and $t_0$ is arbitrary, \eqref{A-2.13} yields
$$
V_1(t_0) =  - C_1 \, t_0 +   \int_\MM \la \nu_{t_0,x}, U \ra \, \theta_1 \, dV_g(x)
\leq
\int_\MM U(u_0) \, \theta_1 \, dV_g(x)
$$
for all $t_0>0$ and, in particular, for all $\theta_1 = \theta_1(x) \geq 0$
\be
\lim_{t \to 0+}  \int_\MM \la \nu_{t,x}, U \ra \, \theta_1 \, dV_g(x)
\leq
\int_\MM U(u_0) \, \theta_1 \, dV_g(x).
\label{A-2.14}
\ee
Note that the left-hand limit exists, in view of (a).

Next, consider the set of all linear, convex and finite combinations
of the form
$$
\sum_j \alpha_j \, \theta_{1,j}(x) \, U_j(u),
$$
where $\alpha_j \geq 0$, $\sum_j \alpha_j = 1$,
the functions $U_j$ are smooth and convex in $u$ and
the functions $\theta_{1,j}(x) \geq 0$ are smooth and compactly support,
with moreover
$$
|U_j(u) \, \theta_{1,j}(x)| \leq c \, |u| + |\widetilde U_j(x)|
$$
with $c \geq 0$ and $\widetilde U_j \in L^1(\MM; dV_g)$.
This set is dense (for the uniform topology in $u$ and the $L^1$ topology in
$x$) in the set of all functions $U=U(u,x)$
that are convex in $u$ and measurable in $x$ and satisfy
$$
|U(u,x)| \leq c \, |u| + |\widetilde U(x)|
$$
for some $c>0$ and $\widetilde U \in L^1(\MM; dV_g)$.
Therefore, by density, (b) follows from \eqref{A-2.13}.

Finally,  (c) follows from (b) by choosing $U(u,x) = |u - u_0(x)|$.
\end{proof}

\begin{proof}[Proof of Theorem~\ref{3-4}]
In all of the following arguments, the inequalities should be understood in the sense
of distributions. All steps can be justified rigorously by introducing test-functions
in the usual way. Restricting attention in \eqref{3.5} to functions $\theta$ with
compact support, we deduce that \be \del_t \la \nu, U \ra + \dive \la \nu, F \ra
\leq 0 \label{3.6} \ee and so, after introducing the Kruzkov's entropies
$$
\Ub(\ubar,\vbar) = |\vbar - \ubar|, \quad \Fb(\ubar,\vbar) = \bigl( f(\vbar) -
f(\ubar) \bigr) \, \sgn(\vbar - \ubar),
$$
we obtain \be \del_t \la \nu, \Ub(\cdot,\vbar) \ra + \dive \la \nu, \Fb(\cdot,\vbar)
\ra \leq 0,  \quad \vbar \in \RR. \label{3.7} \ee Let $\mu$ be another entropy
measure-valued solution to \eqref{2.3}. We are going to combine \eqref{3.7} together
with a similar statement for $\mu$, that is, 
\be 
\del_t \la \mu, \Ub(\ubar,\cdot) \ra + \dive \la \mu, \Fb(\ubar,\cdot) \ra \leq 0,  
\quad \ubar \in \RR. 
\label{3.8} 
\ee

Introducing the tensor product $\nu \otimes \mu = \nu_{t,x} \otimes \mu_{t,x}$ with, for instance,
$$
\la \nu_{t,x} \otimes \mu_{t,x}, \Ub \ra
:= \iint_{\RR^2} \Ub(\ubar,\vbar) \, d\nu_{t,x}(\ubar) d\mu_{t,x}(\vbar),
$$
we can write (in the sense of distributions)
$$
\aligned & \del_t \la \nu \otimes \mu, \Ub \ra + \dive \la \nu \otimes \mu, \Fb \ra
\\
& = \la \nu, \del_t \la \mu, \Ub \ra + \dive \la \mu, \Fb \ra \ra
    + \la \mu, \del_t \la \nu, \Ub \ra + \dive \la \nu, \Fb \ra \ra
\endaligned
$$
and we deduce from \eqref{3.7} and \eqref{3.8} that, in the sense of distributions \be
\del_t \la \nu \otimes \mu, \Ub \ra + \dive \la \nu \otimes \mu, \Fb \ra \leq 0.
\label{3.9} \ee

Next, integrating \eqref{3.9} over the manifold $\MM$ we find that, for all $0 \leq t' \leq t$,
\be
\int_\MM \la \nu_{t,x} \otimes \mu_{t,x}, \Ub \ra  \,  dV_g(x)
\leq
\int_\MM \la \nu_{t',x} \otimes \mu_{t',x}, \Ub \ra \,  dV_g(x).
\label{3.10}
\ee
Letting $t' \to 0$ in \eqref{3.10} and using that the two measure-valued solutions assume
the same (Dirac-mass) initial data $\delta_{u_0(x)}$ at the time $t=0$ in the strong sense established
in Lemma~\ref{3-5}, we find
$$
\int_\MM \la \nu_{t,x} \otimes \mu_{t,x}, \Ub \ra \, dV_g(x) = 0, \quad t \in \RR_+.
$$
Therefore, $\la \nu_{t,x} \otimes \mu_{t,x}, \Ub \ra$ vanishes for almost every $(t,x)$,
which is equivalent to saying that $\nu_{t,x}$ and $\mu_{t,x}$ coincide with the same
Dirac mass, say $\delta_{u(t,x)}$ for some value $u(t,x)$. Since $\nu$ and $\mu$ are
both entropy measure-valued solutions, the function $u=u(t,x)$ is an entropy solution
of the problem \eqref{2.3}-\eqref{2.init} in the sense of Definition~\ref{2-3}.

We establish some additional properties of the entropy solution $u$, as follows.
Taking $U(u) : = |u|^p$ in \eqref{3.6} and integrating over the manifold $\MM$ yields
$$
{d \over dt} \int_\MM U(u(t,x)) \, dV_g(x)
=
{d \over dt} \int_\MM \la \nu_{t,x}, U \ra \, dV_g(x)
\leq  0,
$$
which leads to \eqref{3.1}.
Using \eqref{3.10} with two distinct solutions $\nu_{t,x} = \delta_{u(t,x)}$
and $\mu_{t,x} = \delta_{v(t,x)}$ gives the $L^1$ contraction property \eqref{3.3}.
This completes the proof of Theorem~\ref{3-4}.
\end{proof}


\section{General conservation laws and balance laws}
\label{GE-0}

\subsection*{$L^1$ semi-group of entropy solutions on manifolds}

Consider first the case of a conservation law posed on the one-dimensional torus $T^1= [0,1]$
\be
\del_t u + {1 \over k} \del_x (k \, f(u)) =0, \quad u=u(t,x) \in \RR, \, t \in \RR_+, \, x \in [0,1],
\label{GE.1}
\ee
where $k=k(x)$ is a given, positive function and $f:\RR \to \RR$ is a convex function.
In \cite{LeFlochNedelec}, a suitable generalization of Lax's explicit formula
\cite{Lax} was introduced for \eqref{GE.1}. When $k$ is not a constant, the characteristics of \eqref{GE.1}
are not straight lines but curves $s \mapsto X(s) = X(s; y)$ ($s \geq 0$, $y \in \RR$), given by
$$
\del_s X(s) = \del_u f(u(s, X(s))), \qquad X(0) = y.
$$
It was observed that along a characteristic the function $v=v(s; y) := u(s, X(s; y))$ is such that the ``weighted flux''
$k(X) \, f(v)$ is constant. By introducing suitable left- and right-inverses of the function $f$, say $f^{-1}_\pm$
and then solving the equation $k(X) \, f(v) = c$, it follows that the whole family of all characteristic curves is described by
\be
\del_s X(s) = (\del_u f \circ f^{-1}_\pm)\Big( {c \over k(X)} \Big), \qquad X(0) = y.
\label{GE.curves}
\ee

The following result was derived from a detailled analysis of this family of curves:

\begin{theorem} {\rm (See \cite{LeFloch0,LeFlochNedelec}.)}
\label{LN}
The periodic, entropy solutions with bounded variation to the conservation law \eqref{GE.1}
are given by a generalization of Lax's explicit formula via a minimization problem along the curves \eqref{GE.curves}.
Moreover, any two solutions $u,v$ satisfy the $L^1$ contraction property
$$
\|v(t) - u (t) \|_{L^1(0,1)} \leq  \|v(t') - u (t') \|_{L^1(0,1)}, \qquad 0 \leq t' \leq t.
$$
\end{theorem}

In view of the discussion in the previous sections,
this result may seem surprising since the geometry-compatibility condition is not satisfied here
(except in the trivial case where $k$ is a constant). Theorem~\ref{LN}
motivates us to extend now our theory on Riemannian manifolds
to general conservation laws that need not be geometry-compatible.

From now on we consider a general conservation law \eqref{2.3} associated with an
arbitrary flux $f_x$. First of all, we stress that the notion of entropy pair
$(U,F_x)$ should still be defined by the same conditions \eqref{ent-flux} as in the
geometry-compatible case, but now we no longer have \eqref{ent-sol}. Instead, an entropy
solution should be characterized by the entropy inequalities \be \del_t U(u) + \dive
\big( F(u) \big) - (\dive F)(u) \leq 0 \label{GE.ineq} \ee in the sense of
distributions, for every entropy pair $(U,F_x)$. The term $(\dive F)(\ubar)$
is defined by applying the divergence operator to the vector field $x \to F_x(\ubar)$,
for every fixed $\ubar$.

\begin{theorem}
Let $f=f_x(\ubar)$ be an arbitrary (not necessarily geometry-compatible) flux on a Riemannian manifold $(\MM,g)$,
satisfying the linear growth condition (for some constant $C_0>0$)
\be
| f_x(\ubar) |_g \leq C_0 \, (1 + |\ubar|), \qquad \ubar \in \RR, \, x \in \MM.
\label{growth}
\ee
Then there exists a unique contractive, semi-group of entropy solutions
$$
u_0 \in L^1(\MM) \mapsto u(t) := S_t u_0 \in L^1(\MM)
$$
to the initial value problem \eqref{2.3}-\eqref{2.init}.
\end{theorem}

The condition \eqref{growth} is required for the flux term \eqref{2.3} to be an integrable function on $\MM$.
Note that no uniform $L^p$ estimate is now available,
not even in the $L^1$ norm, since the trivial function $u \equiv 0$
need not be a solution of the conservation law. The $L^1$ norm of a solution is finite for each time $t$,
but generally {\sl grows as $t$ increases.}
This result shows that the contraction property is more fundamental than
all of the other stability properties derived earlier for geometry-compatible conservation laws.

For general conservation laws, the stationary solutions $\widetilde u$,
determined by
$$
\dive \big( f_x(\widetilde u(x))\big) = 0, \qquad x \in \MM,
$$
represent possible asymptotic states of time-dependent solutions.

\begin{proof} The semi-group is constructed first over functions with bounded variation and
then extended by density to the whole of $L^1$. To this end, we re-visit the proof of
Theorem~\ref{3-4}, and check that the key inequality \eqref{3.9},
needed in the derivation of the $L^1$ contraction property, remains true
regardless of the geometry-compatibility property. Note in passing that the stronger statement \eqref{3.6}
is {\sl no longer valid} under the present assumptions, and was precisely needed to ensure the stability
in all $L^p$ spaces, which is no longer true here.

It thus remains to establish  \eqref{3.9} for solutions $u,v$ of \eqref{2.3} with bounded variation.
We rely on standard regularity results (Federer \cite{Federer}, Volpert~\cite{Volpert}) for such functions:
pointwise values $u_\pm, v_\pm$ can be defined almost everywhere with respect to the Hausdorff measure $\Hcal_n$
on $\RR_+ \times \MM$. These pointwise values coincide at points of approximate continuity.
We compute the entropy dissipation measure
$$
\mu := \del_t \Ub(u,v) + \dive\big( \Fb(u,v) \big),
$$
and distinguish between the set $\Ccal_{u,v}$ of points of approximate continuity for both $u$ and $v$,
and the set $\Scal_{u,v}$ of points
where $u$ has an approximate jump and $v$ is approximately continuous, or vice-versa.
The $\Hcal_n$-Hausdorff measure of the set $\Ncal_{u,v} := \RR_+ \times \MM \setminus \Ccal_{u,v} \cup \Scal_{u,v}$ is zero.
Using standard calculus for BV functions we can write
$$
\aligned
\mu_{|\Ccal_{u,v}}
& := \sgn(u-v) \, \big( \del_t u - \del_t v \big) + \sgn(u-v) \, \big( \dive (f(u)) - \dive (f(v)) \big)
\\
& =0.
\endaligned
$$
On the other hand, if $B \subset \Scal_{u,v}$ is a Borel subset consisting of points where
(for instance) $v$ is approximately continuous then
$$
\aligned
&\mu(B) 
\\
&:=  \int_B \Big( n_t \, | u_+ - v| + g_x \Big( n_x, \sgn(u_+ - v) \,
(f_x(u_+) - f_x(v))\Big)
\\
& \quad - n_t \, |u_- - v| - g_x \Big( n_x, \sgn(u_- - v) \, (f_x(u_-) - f_x(v))\Big) \Big) \,  d\Hcal_n(t,x),
\endaligned
$$
where $n=(n_t(t,x), n_x(t,x) \in \RR \times T_x\MM$ is the propagation speed vector
associated with the discontinuity in $u$. Now, by relying on the entropy
inequalities, $\mu\big\{ (t,x) \}$ can easily be checked to be non-positive.
In turn, the measure $\mu$ is also non-positive on $\Scal_{u,v}$ and, in turn,
on the whole of $\RR_+ \times \MM$. The argument is complete.
\end{proof}


\subsection*{General balance laws}
It is not difficult to generalize the well-posed\-ness theory to the balance law
$$
\del_t u + \dive \big( f(u, \cdot) \big) = h(u,\cdot),
$$
where $f=f_x(\ubar,t)$ is a family of vector fields depending (smoothly)
on the time variable $t$ and on the parameter $\ubar$, and $h=h(\ubar,t,x)$ is a smooth function.
Our previous results should be modified to take
into account the dependence of $f$ in $t$ and the effect of the source-term $h$.
For instance, for $L^\infty$ initial data, the $L^p$ estimate \eqref{3.1} should be replaced with
$$
\|u(t)\|_{L^p(\MM; dV_g)} \leq C_p(T) \, \|u(t')\|_{L^p(\MM; dV_g)} + C'(T),
\quad 0 \leq t' \leq t \leq T,
$$
where the constant $C_p(T) \geq 1$ in general depends on $T$
and $C'(T)$ also depends upon $g$. Similarly, the contraction property should be replaced by
$$
\| v(t) - u(t) \|_{L^1(\MM; dV_g)} \leq  C_0(T) \, \| v_0 - u_0 \|_{L^1(\MM; dV_g)}, \qquad t \in \RR_+.
$$
for some constant $C_0(T)$. We omit the details.


\section{Conservation laws on Lorentzian manifolds}
\label{LO-0}

\subsection*{Globally hyperbolic Lorentzian manifold}
Motivated by the mathematical developments in general relativity, we now extend our theory to Lorentzian manifolds.

Let $(\Lorentz,g)$ be a time-oriented, $(n+1)$-dimensional Lorentzian manifold,
$g$ being a pseudo-Riemannian metric tensor on $\Lorentz$ with signature $(-, +, \ldots, +)$.
Recall that tangent vectors $X$ on a Lorentzian manifold can be separated into time-like vectors
($g(X,X) < 0$), null vectors ($g(X,X) = 0$), and space-like vectors ($g(X,X) > 0$). A vector field
is said to be time-like, null, or space-like if the corresponding property hold at every point.
The null cone separates time-like vectors into future-oriented and past-oriented ones; it is assumed
here that the manifold is time-oriented, i.e.~a consistent orientation can be chosen throughout the manifold.

Let $\nabla$ be the Levi-Cevita connection associated with the Lorentzian metric $g$ so that, in particular,
$\nabla g=0$. The divergence $\dive$ operator is defined in a standard way which is formally similar to the Riemannian case.

\begin{definition}
A {\rm flux} on the manifold $\Lorentz$
is a vector field $x \mapsto f_x(\ubar) \in T_x\Lorentz$, depending on a parameter $\ubar \in \RR$.
The {\rm conservation law} on $(\Lorentz,g$) associated with $f$ is
\be
\dive \big( f(u) \big) = 0, \qquad u: \Lorentz \to \RR.
\label{LO.1b}
\ee
It is said to be {\rm geometry compatible} if $f$ satisfies the condition
\be
\dive f_x (\ubar) = 0, \quad \ubar \in \RR, \, x \in \Lorentz.
\label{LO.divergencefree}
\ee
Furthermore, $f$ is said to be a {\rm time-like flux} if
\be
g_x\big(\del_u f_x(\ubar), \del_u f_x(\ubar) \big) < 0, \qquad x \in \Lorentz, \, \ubar \in \RR.
\label{LO.timelike}
\ee
\end{definition}

Note that our terminology here differs from the one in the Riemannian case,
where the conservative variable was singled out.

We are interested in the initial-value problem associated with \eqref{LO.1b}. We fix a space-like hypersurface
$\Hcal_0 \subset \Lorentz$ and a measurable and bounded function $u_0$ defined on $\Hcal_0$.
Then, we search for a function $u=u(x) \in L^\infty(\Lorentz)$ satisfying \eqref{LO.1b}  in the distributional
sense and such that the (weak) trace of $u$ on $\Hcal_0$ coincides with $u_0$:
\be
u_{|\Hcal_0} = u_0.
\label{LO.initial}
\ee
It is natural to require that the vectors $\del_u f_x(\ubar)$, which determine the propagation of waves in
solutions of  \eqref{LO.1b}, are time-like and future-oriented.

We assume that the manifold $\Lorentz$ is {\sl globally hyperbolic,} in the sense that
there exists a foliation of $\Lorentz$ by space-like, compact, oriented hypersurfaces $\Hcal_t$ ($t \in \RR$):
$$
\Lorentz = \bigcup_{t \in \RR} \Hcal_t.
$$
Any hypersurface $\Hcal_{t_0}$ is referred to as a {\sl Cauchy surface} in $\Lorentz$,
while the family $\Hcal_t$ ($t \in \RR$) is called an {\sl admissible foliation associated with} $\Hcal_{t_0}$.
The future of the given hypersurface will be denoted by
$$
\Lorentzplus := \bigcup_{t \geq 0} \Hcal_t.
$$
Finally we denote by $n^t$ the future-oriented, normal vector field to each $\Hcal_t$,
and by $g^t$ the induced metric.
Finally, along $\Hcal_t$, we denote by $X^t$ the normal component of a vector field $X$, thus
$X^t := g(X,n^t)$.

\begin{definition} A flux $F=F_x(\ubar)$ is called a
{\rm convex entropy flux} associated with the conservation law \eqref{LO.1b}
if there exists a convex function $U:\RR \to \RR$ such that
$$
F_x(\ubar) = \int_0^\ubar \del_u U(u') \,\del_u f_x(u') \, du', \qquad x \in \Lorentz, \, \ubar \in \RR.
$$
A measurable and bounded function $u=u(x)$ is called an {\rm entropy solution}
of the geometry-compatible conservation law \eqref{LO.1b}-\eqref{LO.divergencefree}
if the following {\rm entropy inequality}
\be
\int_{\Lorentzplus} g(F(u), \grade \theta) \, dV_g + \int_{\Hcal_0} g_0(F(u_0), n_0) \, \theta_{\Hcal_0} \, dV_{g_0} \geq 0.
\label{LO.3}
\ee
for all convex entropy flux $F=F_x(\ubar)$ and all smooth functions $\theta \geq 0$ compactly supported in
$\Lorentzplus$.
\end{definition}

In particular, in \eqref{LO.3} the inequality
$$
\dive \big( F(u) \big) \leq 0
$$
holds in the distributional sense.

\begin{theorem}
\label{LO-theo}
Consider a geometry-compatible conservation law \eqref{LO.1b}-\eqref{LO.divergencefree}
posed on a globally hyperbolic Lorentzian manifold $\Lorentz$.
Let $\Hcal_0$ be a Cauchy surface in $\Lorentz$, and $u_0: \Hcal_0 \to \RR$ be a measurable and
bounded function. Then, the initial-value problem \eqref{LO.1b}-\eqref{LO.initial}
admits a unique entropy solution $u=u(x) \in L^\infty(\Lorentz)$.
For every admissible foliation $\Hcal_t$, the trace $u_{\Hcal_t}$ exists and belong to $L^1(\Hcal_t)$,
and the functions
$$
\| F^t(u_{\Hcal_t} \|_{L^1(\Hcal_t)},
$$
are non-increasing in time, for any convex entropy flux $F$. Moreover, given any two entropy solutions $u,v$,
the function
\be
\label{LO.contract}
\| f^t(u_{\Hcal_t}) - f^t(v_{|\Hcal_t}) \|_{L^1(\Hcal_t)}
\ee
is non-increasing in time.
\end{theorem}

We emphasize that, in the Lorentzian case, no time-translation property is available in general, contrary to the
Riemannian case. Hence, no time-regularity is implied by the $L^1$ contraction property.

As the proof is very similar to the one in the Riemannian case, we will content with sketching the proof.
Introduce a local chart
$$
x=(x^\alpha)=(t,x^j), \quad g:= g_{\alpha\beta} \, dx^\alpha dx^\beta,
$$
where by convention greek indices describe $0,1, \ldots, n$ and latin indices describe $1,\ldots, n$.
By setting $f =: (f^\alpha_x(\ubar))$ and using local coordinates, the conservation law \eqref{LO.1b} reads
\be
\del_\alpha \big( |g_x|^{1/2} \, f_x^\alpha(u(x)) \big) = 0,
\label{LO.2}
\ee
where $ |g|:= \det (g_{\alpha\beta})$. Thanks to the assumption on $f$, for all {\sl smooth} solutions
\eqref{LO.2} takes the equivalent form
\be
g_x( (\del_u f_x)(u(x)), \grade u(x) \big) := \big(\del_u f_x^\alpha\big) (u(x)) \, (\del_\alpha u)(x) = 0.
\label{LO.22}
\ee

In other words, setting $x= (t,\xbar)$ and $f=(f^t, f^j)$,
$$
\big(\del_u f_x^t\big) (u(t,\xbar)) \, \del_t u(t,\xbar) + \big(\del_u f_x^j\big) (u(t,\xbar)) \, (\del_j u)(t,\xbar) = 0,
$$
in which, since $\del_u f$ is future-oriented and time-like, the coefficient in front of the time-derivative
is positive
$$
\del_u f_x^0(\ubar) > 0, \quad \ubar \in \RR, \, x \in \Lorentz.
$$

To proceed with the construction of the entropy solutions, we add a vanishing diffusion term, as follows:
\be
\dive \big( f(u^\eps) \big) = \eps \, \Deltabar_{\overline{g}} u^\eps,
\label{LO.diff}
\ee
where $\Deltabar_{\overline{g}}$ is the Laplace operator on the leaves $\Hcal_t$ of the foliation, that is
in coordinates
\be
\label{LO.local}
\aligned
&\big(\del_u f_x^0\big) (u^\eps(t,\xbar)) \, \del_t u^\eps(t,\xbar)
   + \big(\del_u f_x^j\big) (u^\eps(t,\xbar)) \, (\del_j u^\eps)(t,\xbar)
\\
& = \eps \, g^{ij}(t,\xbar) \, \big( \del_i \del_j u^\eps - \Gamma^k_{ij} \, \del_k u^\eps \big)(t,\xbar).
\endaligned
\ee

In view of \eqref{LO.local}, we see that all of the estimates follow similarly as in the Riemannian case.
For instance, multiplying \eqref{LO.local} by a convex function $U$, the entropy inequality now takes the form
\be
\label{LO.entrop}
\aligned
& \del_t (F_x^0 (u^\eps(t,\xbar)) + \del_j (F_x^j (u^\eps(t,\xbar))
\\
& = \del_u U (u(t,\xbar)) \, \big(\del_u f_x^0\big) (u^\eps(t,\xbar)) \, \del_t u^\eps(t,\xbar)
\\
& \qquad   + \del_u U (u(t,\xbar)) \, \big(\del_u f_x^j\big) (u^\eps(t,\xbar)) \, (\del_j u^\eps)(t,\xbar)
 \\
& = \eps \, \del_u U (u(t,\xbar)) \, g^{ij}(t,\xbar) \, \big( \del_i \del_j u^\eps - \Gamma^k_{ij} \, \del_k u^\eps \big)(t,\xbar)
\\
& = \eps \, \Deltabar_{\overline{g}} U(u^\eps) - \eps \, \del_u^2 U(u^\eps) \, g^{jk} \, \del_j u^\eps \del_k u^\eps.
\endaligned
\ee
The metric $\overline g$ induced on the space-like leaves $\Hcal_t$ is positive-definite and, therefore,
the latter term above if non-positive and we conclude that, given a geometry compatible flux $f$
and for every convex entropy flux $F$
$$
\dive \big( F(u^\eps) \big)
=
\del_t (F_x^0 (u^\eps(t,\xbar)) + \del_j (F_x^j (u^\eps(t,\xbar))
\leq
\eps \, \Deltabar_{\overline{g}} U(u^\eps).
$$

Note that the regularization \eqref{LO.diff}
does depend on the specific foliation under consideration. However, by the contraction property
\eqref{LO.contract} the limiting solution is unique and independent of the chosen regularization
mechanism.


\subsection*{Schwarzschild spacetime}

In the context of general relativity, the Schwar\-schild metric represents a spherically symmetric empty
spacetime surrounding a black hole with mass $m$ and is one of the most important example of Lorentzian metrics.
The outer communication region of the Schwarz\-schild spacetime
is a $1+3$-Lorentzian manifold with boundary, described
in the so-called Schwarschild coordinates $(t,r,\omega)$ by
$$
g = -\bigl( 1 - {2m \over r} \bigr) \, dt^2
  + \bigl( 1 - {2m \over r} \bigr)^{-1} \, dr^2 + r^2 \, d\omega^2,
$$
with $t>0$ and $r > 2 m$, while $\omega$ describes the $2$-sphere.
There is an apparent (but not a physical) singularity in the metric coefficients at $r= 2m$,
which corresponds to the horizon of the spacetime.
This spacetime is spherically symmetric, that is invariant under the group of
rotations operating on the space-like $2$-spheres given by keeping $t$ and $r$ constant.
The part $r^2 d\omega^2$ of the metric is the canonical metric on the $2$-spheres of symmetry.
The spacetime under consideration is static, since the vector field $\del_t$ is a time-like Killing vector.
Moreover, this metric is asymptotic to the flat metric when $r \to \infty$.
Theorem~\ref{LO-theo} extends to the exterior of the Schwarzschild spacetime, by
observing that along the boundary $r=2m$ the characteristics of the hyperbolic equation are outgoing.


\section*{Acknowledgments}
After this paper was submitted, the authors learned about an earlier paper on the
subject by E.Y. Panov, On the Cauchy problem for a first-order quasilinear equation on a manifold,  
Differential Equations 33 (1997), 257--266. 

The authors were supported by a
research grant ``Arc-en-Ciel'' sponsored by the High Council for Scientific and
Technological Cooperation between France and Israel,  entitled: {\em ``Theoretical and
numerical study of geophysical fluid dynamics in general geometry''.} P.G.L. was also
supported
by the A.N.R. grant 06-2-134423 entitled {\sl ``Mathematical Methods in General Relativity''} (MATH-GR).


\newcommand{\auth}{\textsc}


\begin{thebibliography}{10}

\bibitem{AmorimBenArtziLeFloch} \auth{P. Amorim, M. Ben-Artzi, and P.G. LeFloch},
Hyperbolic conservation laws on manifolds.
Total variation estimates and the finite volume method,
{\em Meth. Appli. Anal.} {\bf 12} (2005), 291--324.

\bibitem{Dafermos} \auth{C.M. Dafermos,}
{\sl Hyperbolic conservation laws in continuum physics,}
Grundlehren Math. Wissenschaften Series, Vol. {\bf 325}, Springer Verlag, 2000.

\bibitem{DiPerna} \auth{R.J. DiPerna,}
Measure-valued solutions to conservation laws,
{\em Arch. Rational Mech. Anal.} {\bf 88} (1985), 223--270.

\bibitem{Federer} \auth{H. Federer,}
{\sl Geometric measure theory,} Springer-Verlag, New York, 1969.

\bibitem{Hormander} \auth{L. H\"{o}rmander,}
{\sl Nonlinear hyperbolic differential equations,} Math. and Appl. {\bf 26,} Springer Verlag, 1997.

\bibitem{Kruzkov} \auth{S. Kruzkov,}
First-order quasilinear equations with several space variables,
{\em Math. USSR Sb.} {\bf 10} (1970), 217--243.

\bibitem{Lax} \auth{P.D. Lax,}
{\sl Hyperbolic systems of conservation laws and the mathematical theory of shock waves,}
Regional Conf. Series in Appl. Math. {\bf 11}, SIAM, Philadelphia, 1973.

\bibitem{LeFloch0} \auth{P.G. LeFloch},
Explicit formula for scalar conservation laws with boundary condition,
{\em Math. Meth. Appl. Sc.} {\bf 10} (1988), 265--287.

\bibitem{LeFloch} \auth{P.G. LeFloch},
{\sl Hyperbolic systems of conservation laws: The theory of classical
and nonclassical shock waves},
Lectures in Mathematics, ETH Z\"urich, Birkh\"auser, 2002.

\bibitem{LeFlochNedelec} \auth{P.G. LeFloch and J.-C. Nedelec,}
Explicit formula for weighted scalar nonlinear conservation laws,
{\em Trans. Amer. Math. Soc.} {\bf 308} (1988), 667--683.

\bibitem{Spivak} \auth{M. Spivak},
{\sl A comprehensive introduction to differential geometry,}
Vol.~{\bf 4}, Publish or Perish Inc., Houston, 1979.

\bibitem{Volpert} \auth{A.I. Volpert,} The space BV and quasi-linear equations,
{\em Mat. USSR Sb.} {\bf 2} (1967), 225--267.

\end{thebibliography}
\end{document}